%% file: artsum3.tex
\documentclass[a4paper]{article}
\input{headart.tex}

\begin{document}
\title{Approximation of *weak-to-norm continuous mappings}
\author{Lorenzo D'Ambrosio\thanks{SISSA-ISAS,
via Beirut 2--4,  34014  Trieste Italy.\ \ 
e-mail address: dambros@sissa.it}}
\date{\fs Ref. SISSA 75/2000/M}
\maketitle\label{ver 3.7}
\pagenumbering{arabic}
\input{1sec3}

\input{2sec3}
\input{3sec3}
\input{4sec4}
\bibliography{bibjor,bibbook}
\bibliographystyle{lortesi}
\end{document}

%% file: headart.tex
\usepackage{latexsym}

\newcommand{\be}{\begin{equation}}	%\be  :=\begin equation
\newcommand{\ee}{\end{equation}}	%\ee  :=\end equation
\newcommand{\bern}{\begin{eqnarray*}}	%\berm  :=\begin equationarray no number
\newcommand{\eern}{\end{eqnarray*}}	%\eern  :=\end equationarray no number
\newcommand{\beqp}{\begin{eqproof}}	%\beqp  :=\begin equation proof
\newcommand{\eeqp}{\end{eqproof}}	%\eeqp  :=\end equation proof

\newcommand{\bt}{\begin{teorema}}	%\bt  :=\begin teorema
\newcommand{\et}{\end{teorema}} 	%\et  :=\end teorema
\newcommand{\bd}{\begin{definizione}}	%\bd  :=\begin definizione
\newcommand{\ed}{\end{definizione}}	%\ed  :=\end definizione
\newcommand{\bc}{\begin{corollario}}	%\bc  :=\begin corollario
\newcommand{\ec}{\end{corollario}}	%\ec  :=\end corollario
\newcommand{\bp}{\begin{dimostrazione}}	%\bp  :=\begin dimostrazione
\newcommand{\ep}{\end{dimostrazione}}	%\ep  :=\end dimostrazione
\newcommand{\bl}{\begin{lemma}}		%\bp  :=\begin lemma
\newcommand{\el}{\end{lemma}}		%\ep  :=\end lemma
\newcommand{\bpr}{\begin{proposizione}}	%\bp  :=\begin proposizione
\newcommand{\epr}{\end{proposizione}}	%\ep  :=\end proposizione
\newcommand{\besi}[1][{}]
	{\begin{esempi}\rm\textbf{#1}\begin{enumerate}}	%\besi  :=\begin esempi
\newcommand{\eesi}{\end{enumerate}\end{esempi}}	%\eesi  :=\end esempi
\newcommand{\beso[1]}{\begin{esempio}[#1]\rm}	%\beso  :=\begin esempio
\newcommand{\eeso}{\end{esempio}}	%\eeso  :=\end esempio
\newcommand{\boss}{\begin{osservazione}\rm}	%\boss  :=\begin osservazione
\newcommand{\eoss}{\end{osservazione}}	%\eoss  :=\end osservazione
\newcommand{\refe}[1]{(\ref{#1})}	%\refe	:=riferimento tra parantesi

\font\corsivo=rsfs10 %at 12pt
\font\doppio=msbm10 %at 12pt

\font\fs=cmr9

\newcommand{\ind}{\!\!}
\newcommand{\alli}{\ind &\ind}

\newcommand{\B}{\mbox{\corsivo B}}	%\B	:=simb dei Boreliani
\newcommand{\C}{\mbox{\corsivo C}}	%\C	:=simb di continuita
\newcommand{\F}{\mbox{\corsivo F}}	%\F	:=simb di funzioni
\newcommand{\K}{\mbox{\corsivo K}}	%\K	:=simb di spazi K
\newcommand{\lin}{{\mathcal L}}		%\lin	:=simb di spazio di linearit
\newcommand{\R}{\hbox{\doppio R}}	%\R	:=simb dei Reali
\newcommand{\uni}{\mathbf{1}}		%\uni	:=simb della funz costante 1

\newcommand{\lb}{[\![}			%\lb  :=doppia quadra
\newcommand{\rb}{]\!]}			%\rb  :=doppia quadra
\newcommand{\norm}[1]{\left\|#1\right\|}		%\norm{}  :=norma di{}
\newcommand{\decl}{:=}			%\decl  :=uguale per definizione
\newcommand{\de}{\mathrm{d}}
\newcommand{\risp}[1]{\textrm{\Big[resp.}#1\textrm{\Big]}}

\renewcommand{\theequation}{\thesection.\arabic{equation}}

\newtheorem{teorema}{Theorem}[section]
\newtheorem{definizione}[teorema]{Definition}
\newtheorem{corollario}[teorema]{Corollary}
\newtheorem{lemma}[teorema]{Lemma}
\newtheorem{proposizione}[teorema]{Proposition}
\newtheorem{esempio}[teorema]{Example}
\newtheorem{esempi}[teorema]{Examples}
\newtheorem{osservazione}[teorema]{Remark}

\newcounter{proofeqno}		%contatore di equazione nelle dimostrazioni
\newenvironment{dimostrazione}	%ambiente per la dimostrazioni
	{\noindent \textbf{Proof.\ }\setcounter{proofeqno}{0}}
	{\hfill $\Box$\par\medskip}

\newenvironment{eqproof}	%come equation, ma nelle dimostrazioni
	{\addtocounter{proofeqno}{1}$$}
	{\eqno(\theproofeqno)$$}

	{\addtocounter{proofeqno}{1} \addtocounter{equation}{-1}
		\renewcommand\theequation{\theproofeqno}
		\begin{eqnarray}   }
	{ \end{eqnarray}
		\renewcommand\theequation{\thesection.\@arabic\c@equation}	}

	{\mbox{}\par }{\par\medskip}

\makeindex

%% file: 1sec3.tex
%\title{\bf Approximation of vector-valued mappings}
%\author{Lorenzo D'Ambrosio}
%{\thanks
%\date{}
%\maketitle
%%%%%%%%%%%%%%%%%%%%%%%%%%%%%%%%%%%%%% version 2.3
%{\centering {\scriptsize SISSA-ISAS
%
%via Beirut 2--4,  34014  Trieste Italy.
%
%e-mail address: dambros@sissa.it}}
%\vskip1cm

{\abstract{The purpose of this paper is to study the approximation of vector
valued mappings defined on a subset of a normed space.
We investigate Korovkin-type conditions
under which a given sequence of linear operators becomes a so-called
approximation process.
First, we give a sufficient condition for this sequence to approximate the class
of bounded, uniformly continuous functions. Then we present some sufficient
and necessary conditions guaranteeing the approximation within the class
of unbounded, *weak-to-norm continuous mappings.
We also derive some estimates of the rate of convergence.

%{In this paper we deal with the sequences of linear operators which
%under the Korovkin type conditions became the so-called approximation
%processes. The operators act on domains consisting of completely continuous
%functions between Banach spaces.% (possibly of infinite dimension).
%
%%In the first part of the article
%We present a number of conditions, sufficient
%and necessary for a given operators' sequence to be an
%approximation process. We also derivate estimates of convergence rates.
%% (in finite dimensional case).

}}

%%In the last two sections we show some possible applications
%%of the previous results; in particular, making use of the estimates
%%mentioned above, we establish new representation formulae for
%%semigroups of operators.}}

\section{Introduction}
Korovkin's well known result, \cite{korovkin}, states that if $(L_n)_{n\ge 1}$
is a sequence of positive linear operators on $\C(\lb a, b\rb)$ then
$\norm{L_n(f)-f}_\infty\rightarrow 0$ for every $f\in\C(\lb a,b\rb)$,
provided the same is true for the following \emph{test functions}: $f(u)=1$,
$f(u)=u$, $f(u)=u^2$.
Shisha and Mond in \cite{shisha-mond} present a quantitative version of
Korovkin's theorem, containing some estimates of the rate of convergence of
$\norm{L_n(f)-f}$ in terms of the corresponding rate of convergence computed
for the test functions.
Many authors have contributed to understanding the possible enlargement
of the domain of approximation operators, in particular to include classes
of unbounded functions.
%% and some quantitative results have been found.
Ditzian in \cite{ditzian} deals with continuous real valued functions,
defined on a closed and unbounded subset of the real line,
which satisfy the growth condition $|f(u)|\le M_f(1+u^2)\mu(u)$ with $\mu\ge1$.
He estimates the rate of the approximation in terms of the rate of convergence for the
test functions $1$,  $u$, $u^2$ and $(u-t)^2\mu(u)$.
Shaw and Yeh in \cite{shaw-yeh} study the case of functions
defined on an open interval $\rb a,b\lb$ of $\R$ and satisfying
$|f(u)|=O(g_a(u)) (u\rightarrow a^+)$ and $|f(u)|=O(g_b(u)) (u\rightarrow b^-)$
(for some suitable convex functions $g_a$ and $g_b$).
The test functions determining the convergence rates are now the following:
$1$, $u$, $u^2$, $g_a$ and $g_b$.
Shaw in \cite{shaw} considers continuous functions on $\R^m$
with a prescribed growth at infinity. More precisely, he treats operators
$L_n$ defined by means of measures: $L_n(f)(t)=\int f(u)\de\mu_{n,t}(u)$,
and the following classes of functions $f$. The first class consists of those
real valued functions whose growth is controlled
by a convex function $g$. The second admissible class contains
functions of the form $T(u)x$, where $x$ belongs to a Banach space $E$,
and $T(u)$ is a linear continuous operator from $E$ into itself such that
$T(u)$ is bounded on bounded subsets of $\R^m$ and $\norm{T(u)}\le Mg(u)$.
Many authors have also studied the case of vector-valued mappings defined on a
compact Hausdorff space $X$ see e.g.
\cite{nishi92, nishi96} and \cite{prolla}.
The former studies the convergence of a net of quasi-positive linear
operators to an operator $T$, that can be the identity on $\C(X;E)$.
Actually, in \cite{nishi96} the value space $E$ is a Dedekind
complete normed vector lattice with normal unit order and, in \cite{nishi92}
$E$ is a normed linear space. Always in the setting of compactness of $X$,
Prolla studies the approximation processes for the identity on $\C(X;E)$ by
monotonically regular operators (that is the operators that are $S$-regular with $S$
positive, see section 2). Moreover, he gives a rate of approximation
when $X$ is a compact subset of a normed space and the process is made of
dominated operators.

The purpose of this article is to give a generalization of the above results
for classes of mappings defined on a convex subset of a vector space taking
their values into a normed space.
The paper is organized as follows.

In section 2 we introduce the notation and definitions used in the sequel.

Replacing the previous assumption on the positivity of the operators $L_n$ by
the concept of so-called dominated operators we proceed to find Korovkin-type
conditions, as described in the third section.
We also derive there a Korovkin type theorem on the approximation process
within the class of
bounded and uniformly continuous functions defined on a convex set, and
find an estimate of the rate of convergence. In the end of the section
we deduce a Korovkin-type theorem for *weak-to-norm continuous maps on
bounded sets.

The last section deals with the case of unbounded functions.
With $X$ being a *weakly closed or open convex subset of a dual space $Y=Z'$,
we present a Korovkin type theorem for *weak-to-norm
continuous maps on $X$, whose growth is controlled by a convex function.
Under the additional assumption of the dimension of $Y$ to be finite,
we establish some new estimates of the rate of convergence.

%% file: 2sec3.tex
\section{Notation and preliminary definitions}%%%%%%%%ver 2.7
\label{sec2}

In this work $Y$, $Z$ and $E$ will denote real or complex normed spaces,
with their norms denoted by same symbol
$\norm\cdot$. As usual, $Z'$ stands for the dual space of $Z$ and
$\pi(Z)$ stands for the dual space of $Z'$ with *weak topology $\sigma(Z',Z)$,
so $Z$ is reflexive if and only if $\pi(Z)=Z''$.
If $\phi\in \pi(Z)$, and $X$ is a nonempty subset of $Z'$, then by $\phi_{|X}$
we mean the restriction of $\phi$ to $X$.

We will often address to the following two functional spaces:
$\F(X;E)$ and $\B(X;E)$ that are, respectively, the vector space of all mappings
$F\colon X\to E$ and its subspace containing only the bounded mappings.
The latter space is normed by the uniform norm $\norm{\cdot}_X$
%defined, for all \(F\in\B(X;E)\), as
\[\norm{F}_{X}\decl\sup_{u\in X}\norm{F(u)}.\]
For $F$ belonging to the former space, $\norm{F}\colon X\to \R$ denotes
the real valued function
$\norm{F}(u)\decl\norm{F(u)}$.

With the usual symbol $\C(X;E)$ we denote the subspace of $\F(X;E)$ consisting
of all continuous mappings.

Fix $g\colon X\to\R$ a strictly positive function. Then $\C(X;E,g)$ denotes
the subspace of all mappings $F\in\C(X;E)$ such that
$\norm{F(u)}\le Mg(u)$ for every $u\in X$ and some constant $M>0$,
depending only on $F$.
Finally, $UCB(X;E)$ is the subspace of all mappings of  $\C(X;E)$ which are
uniformly continuous and bounded.

In case $E=\R$ we abbreviate the above notation, writing
$\C(X,g)$ instead of $\C(X;\R,g)$, $\F(X)$ instead of $\F(X;\R)$ and so on.

We also adopt the following notation: if $c\in E$, then, we shall denote
again by $c$ the constant mapping $F(u)=c$ ($u\in X$).

If $f\in\F(X)$ and $x\in E$, $f\otimes x$ denotes the mapping of $\F(X;E)$
defined by $(f\otimes x)(u)\decl f(u)x$ $(u\in X)$.

For $t\in Y$, define $\psi_t\colon X\to \R$ by the formula
$\psi_t(u)\decl\norm{u-t}$.
Observe that if $\psi_{t_0}^2\in\C(X,g)$, for some $t_0\in Y$,
then the same holds for every $t\in Y$.

\bd Let $Z$ be normed space, $Y$ its dual space and $X\subset Y=Z'$.
	We say that $F\colon X\to E$ is
	\emph{*weak-to-norm continuous} if it is continuous
	from $X$ equipped with the *weak topology $\sigma(Y,Z)$ in $Y$,
	into $E$ with the norm topology.
	By $\K(X;E)$ we denote the space of all *weak-to-norm continuous
	mappings from $X$ into $E$.
	We set $\K(X;E,g)\decl \K(X;E)\cap \C(X;E,g)$.
\ed

We remark that every *weak-to-norm continuous mapping is in particular
continuous and maps *weakly closed and bounded subsets of $X$
in compact subsets of $E$.
Moreover, if the dimension of $Y$ is finite, then obviously $\K(X;E)=\C(X;E)$.

%If $X$ is convex
For $F\in UCB(X;E)$, as usual, we denote with
$\omega(F,\cdot)$ its modulus of continuity,
%% for all $h\geq0$,
\[\omega(F,h)\decl\sup\{\norm{F(u)-F(t)}\big|\, t,u\in X,\, \norm{t-u}\le h \}\ \  (h>0).\]

The following definitions are based on the analogous ones in \cite{prolla}.

\bd\label{def:oper.dominato} $\!\!$ Let $L\colon D(L)\to\F(X;E)$
	and $S\colon D(S)\to\F(X)$ be linear operators defined on some subspace
	$D(L)$ and $D(S)$ of $\C(X;E)$ and $\C(X)$, respectively.
	We say that
\begin{enumerate}
\item[a)] $L$ is \emph{dominated} by $S$ if
	$\norm F \in D(S)$, and 
	\[\norm{L(F)(t)}\le S(\norm{F})(t)\]
	for all $F\in D(L)$ and $t\in X$;
\item[b)] $L$ is \emph{$S$-regular} if
	$f\otimes x\in D(L)$ and
	\[L(f\otimes x) = S(f)\otimes x \]
	for all $f\in D(S)$ and $x\in E$;

\item[c)]  $L$ \emph{preserves the constants} if
	$c\in D(L)$ and $L(c)(t)=c$, for all $c\in E$ and $t\in X$.
\end{enumerate}
\ed

Below we present some examples of dominated and regular operators.

\beso[Interpolation Operators] Let $\lin(E)$ be the Banach algebra of the continuous
	linear operators on $E$ and $I$ be an index set.
	For every $i\in I$ fix a point $t_i\in X$ and an application
	$\Phi_i\in\C(X;\lin(E))$, and set
	$\phi_i\decl\norm{\Phi_i}_{\lin(E)}\in\C(X)$.
	We consider the operators $L\colon D(L)\to\F(X;E)$ and
	$S\colon D(S)\to\F(X)$, defined by
%%	 come\footnote{Se l'insieme degli indici 
%%	non \`e finito allora la somma \`e da intendersi nel seguente modo. 
%%	Siano $x_i\in E$ e $x\in E$, scriveremo $\sum_{i\in I}x_i=x$ se, e solo se, per
%%	ogni $\epsilon>0$ esiste $I_0\subset I$ finito tale che per ogni $J$ finito
%%	e $I_0\subset J\subset I$ risulta $\norm{\sum_{i\in J}x_i-x}\le\epsilon$.
%%
%%	Con tale definizione le propriet\`a algebriche della somma vengono mantenute.
%%	Inoltre \`e facile verificare, attraverso la definizione, che sussiste la disuguaglianza
%%	\[\norm{\sum_{i\in I}x_i}\le\sum_{i\in I}\norm{x_i}.\]}% Fine FootNote
%%
%%
	\bern L(F)(t)\decl\sum_{i\in I}\langle\Phi_i(t),F(t_i)\rangle&&
					\mathrm{for\ any\ }F\in D(L),\\
	S(f)(t)\decl\sum_{i\in I}\phi_i(t)f(t_i)&&\mathrm{for\ any\ }f\in D(S),\eern
	for all $t\in X$. The domain $D(S)$ is the space of those functions
	$f\in\C(X)$ for which the family $(\phi_i(t)f(t_i))_{i\in I}$ is summable
	for all $t\in X$.
	The domain $D(L)$ is the space of the maps $F\in\C(X;E)$ such that
	$\norm F \in D(S)$. The inequality
	\[\norm{L(F)(t)}\le\sum_{i\in I}\norm{\Phi_i(t)}_{\lin(E)}\norm{F(t_i)}=
			S(\norm{F})(t)\]
	implies that $L$ is well defined on $D(L)$ and that $L$
	is dominated by $S$.

	If for every $i\in I$,  there exists $\psi_i\in \C(X)$ such that
	$\Phi_i(t)(v)=\psi_i(t)v$ ($t\in X, v\in E$),
	then setting $\phi_i\decl\psi_i$, we have
	that $L$ is $S$-regular. Moreover, if $\psi_i\ge 0$ then
	$L$ is also dominated by $S$.
%%	Quindi, da $\phi_i=\norm{\Phi_i}=\psi_i$, ho le relazioni
%%	\bern L(F)(t)\decl\sum_{i\in I}\psi_i(t)F(t_i)&&
%%				\mathrm{per\ ogni\ }F\in D(L);\\
%%	S(g)(t)\decl\sum_{i\in I}\psi_i(t)g(t_i)&&\mathrm{per\ ogni\ }g\in D(S).\eern
\eeso
\beso[Integral Operators]\label{ese:oper.dominato.integr} Let $(E,\norm{\cdot})$
	be a Banach space and assume that for any $t\in X$, a positive
	finite measure $\mu_t\colon\B_X\to\R_+$ on the $\sigma$-~algebra of
	all Borel subset of $X$ is given.
	Define $D(L)\decl\C(X;E)\cap \bigcap_{t\in X}L^1(\mu_t;E)$, and
%	analogously
 $D(S)\decl\C(X)\cap \bigcap_{t\in X}L^1(\mu_t)$.
	Consider the operators  $L\colon D(L)\to\F(X;E)$ and $S\colon D(S)\to\F(X)$
	given by
	\bern L(F)(t)\decl\int_X F(u)\de\mu_t(u)&&\mathrm{for\ any\ }F\in D(L),\\
	S(g)(t)\decl\int_X g(u)\de\mu_t(u)&&\mathrm{for\ any\ }g\in D(S),\eern
	for all $t\in X$.
	Trivially, $L$ and $S$ are linear and $S$ is positive.

	$L$ is dominated \emph{in natural way} by $S$:
	\[\norm{L(F)(t)}=\norm{\int_X F(u)\de\mu_t(u)}\le
		\int_X\norm{F(u)}\de\mu_t(u)=S(\norm{F})(t).\]

	Using the above estimate, we note that for an arbitrary $F\in\C(X;E)$,
	$S(\norm F)$ is well-defined provided $L(F)$ is defined.

	By properties of the Bochner integral it is easy to verify that $L$ is $S$-regular.

	Moreover, we observe that $L$ preserves the constants if and only if
	the measures $\mu_t$ have unit masses or, equivalently,
	$S(\uni)(t)=1$ for all $t\in X$.
\eeso

We will also make use of the following notation:
if $\psi_t^2\in D(S)$ then we write $\gamma^2(t)\decl S(\psi_t^2)(t)$.

%% file: 3sec3.tex
\section{A Korovkin-type theorem for bounded uniformly continuous mappings
between normed spaces}%%%%%%%%%%%%%%%%%% ver 2.4 -----------
\label{sec3}

In this section we approximate vector valued, bounded and uniformly
continuous mappings defined on a convex subset of a normed space.

\bt\label{teo:shimond.vett} Let $Y$ and $E$ be normed spaces,
	$X$ a convex subset of $Y$ and $L_n\colon D(L_n)\to\F(X;E)$ a
	sequence of linear operator dominated by some positive linear operators
	$S_n\colon D(S_n)\to\F(X)$. We suppose that, for every $n\ge1$
%	come nella definizione \ref{def:oper.dominato}	supponendo inoltre che
	$UCB(X;E)\subset D(L_n)$, $UCB(X)\subset D(S_n)$
	and $\psi_t^2\in D(S_n)$ for some (and hence for all) $t\in Y$.
	Then for each $F\in UCB(X;E)$, $t\in X$ and $\delta>0$ one has
	\begin{eqnarray}\lefteqn{\!\!\!\!\!\!\!\!\!\!\!\!\! \norm{L_n(F)(t)-F(t)}\le\norm{L_n(F(t))(t)-F(t)}+
							S_n(\norm{F-F(t)})(t)\nonumber}\\
	&&\qquad\quad\le\norm{L_n(F(t))(t)-F(t)}+\omega(F,\delta)\left[S_n(\uni)(t)+
		\delta^{-2}\gamma_n^2(t)\right]\!\label{dis:shimond.vett} \end{eqnarray}
	where $\gamma_n^2(t)\decl S_n(\psi_t^2)(t)$.

	Moreover if $L_n$ preserves the constants, then:
	\[\norm{L_n(F)(t)-F(t)}\le\omega(F,\delta)\left[S_n(\uni)(t)+
		\delta^{-2}\gamma_n^2(t)\right].\]
	In particular, taking $\delta=\gamma_n(t)$ we obtain
	\[\norm{L_n(F)(t)-F(t)}\le\omega(F,\gamma_n(t))\left[S_n(\uni)(t)+1\right],\]
	and if $\gamma_n$ and $S_n(\uni)$ are bounded on $K\subset X$, then:
	\[\norm{L_n(F)-F}_K\le\omega(F,\norm{\gamma_n}_K)\left[\norm{S_n(\uni)}_K+1\right].\]
\et
\bp Fix $F\in UCB(X;E)$. For every $u\in X$ and $\delta>0$, by the definition
	of $\omega(F,\cdot)$, we get the inequality:
	\[\norm{F(u)-F(t)}\leq\omega(F,\norm{t-u})\leq(1+
		\delta^{-2}\norm{u-t}^2)\,\omega(F,\delta).\]
	Applying the positive operator $S_n$ we have:
	\[ S_n(\norm{F-F(t)})(t)\le\omega(F,\delta)\left(S_n(\uni)(t)+
		\delta^{-2}\gamma_n^2(t)\right),\]
	and
	\bern \norm{L_n(F)(t)-F(t)}&\le&\norm{L_n(F-F(t))(t)}+\norm{L_n(F(t))(t)-F(t)}\\
	&\le& S_n(\norm{F-F(t)})(t)+\norm{L_n(F(t))(t)-F(t)},\eern
	as  $L_n$ is dominated by $S_n$.
\ep

%Applicando il teorema \ref{teo:shimond.vett} potremmo ottenere processi di
%approssimazione anche attraverso successioni di operatori $L_n$ che agiscono su spazi di
%funzioni definite su insiemi non necessariamente limitati.
%Attraverso il teorema
Note that the Theorem \ref{teo:shimond.vett} yields the uniform convergence
of $(L_n(F))_{n\ge1}$ to $F$ on those subsets of $Y$ where the sequence
$\gamma_n^2(t)=S_n(\psi^2_t)(t)$ converges to $0$ uniformly.

When $X$ is *weakly closed, convex and bounded subset of the dual space $Y=Z'$,
	then by Theorem \ref{teo:shimond.vett} and the inclusions
	\bern \K(X;E)\subset UCB(X;E), \\
	\psi^2_t\in UCB(X)\subset D(S_n),\eern
	one obtain the following Korovkin-type theorem for *weak-to-norm continuous maps.
%%\footnote{Proof of $\K(X;E)\subset UCB(X;E)$.
%%	By $I_{\delta,\ell}$ we indicate a neighborhood of the $0$ in the *weak topology
%%	(see Lemma \ref{lem:stime}). First we define the space of uniformly *weak continuous mappings: 
%%	\[w-UCB(X;E)\decl\{F\colon X\to E| \forall\epsilon>0 \ \exists I_{\delta,\ell}
%%	 \ni' x-y\in I_{\delta,\ell}:\norm{F(x)-F(y)}\le\epsilon \}.\] }

\bc\label{cor} Let $Z$ and $E$ be normed spaces, $Y$ the dual space of $Z$,
	$X$ a *weakly closed, convex and bounded subset of the dual space $Y=Z'$,
	and $L_n\colon D(L_n)\to\F(X;E)$ a
	sequence of linear operator dominated by some positive linear operators
	$S_n\colon D(S_n)\to\F(X)$. We suppose that, for every $n\ge1$
	$UCB(X;E)\subset D(L_n)$, $UCB(X)\subset D(S_n)$ and set
	$\gamma_n^2(t)\decl S_n(\psi_t^2)(t)$.
	If for every $c\in E$ the following convergences hold
	\bern L_n(c)\rightarrow c \qquad\mathrm{[resp.\ uniformly\ in\ c\in E];}\\
		\gamma_n(t)\rightarrow 0 \qquad\mathrm{[resp.\ uniformly\ in\ t\in X];}\eern
	then for each $F\in \K(X;E)$
	\[ L_n(F)(t)\rightarrow F(t)\qquad \mathrm{[resp.\ uniformly\ on\ X]}\]
	and moreover the inequalities of the Theorem \ref{teo:shimond.vett} hold.
\ec

\boss In the setting of Corollary \ref{cor}, $X$ results to be
a compact space with *weak topology and, in order to study the approximation process
of the identity on $\K(X;E)$, the above result is slightly different from
the analogue  in \cite[Theorem 1]{prolla} and in
\cite[Corollary 5 and Remark 4]{nishi92}.
Prolla, dealing with dominated operators, requires that $(X,d)$ is a
metric space and the test functions depend on the metric $d$.
In our case, of *weak-to-norm continuous mappings, this means to require
the separability of $Z$ and to use the metric $d$,
%%generated by the numerable dense set of $Z$
given for every $x, y \in X$ by:
\[d(x,y)\decl\sum_{n\ge 1} \frac{|\langle x-y,f_n\rangle|}{2^n},\]
where $f_n\in Z, \norm{f_n}=1$ and $(f_n)_{n\ge1}$ is dense on the unitary sphere
of $Z$. In Corollary \ref{cor} one does not need the separability
of $Z$, and the test functions are based on the easier to use norm of the space.
Nishishiraho tests the sequences of quasi-positive operators on a greater test set
that in our context is
\[\{c \phi_{|X}^k\, |\, \phi\in \pi(Z),\ k=0,1,2\ \mathrm{and }\ c\in E\}.\]
\eoss

The cases of $X$ closed and unbounded, or open are treated in the next section.

%% file: 4sec4.tex
%%%%%%%%%%%%%%%%%%%%%%%%%%%%%%%%% ver 2.5 -----------
\section{Korovkin-type theorems for unbounded mappings between normed spaces}\label{sec4}

As in the scalar case, where it is necessary to control the growth of the
approximated functions
(cfr. \cite{ditzian}), for vector-valued  mappings defined on subsets of
Banach spaces we will have to assume appropriate conditions estimating
the growth near the boundary of their domains of definition.

Since now we assume that $(Z,\norm\cdot)$ is a real normed space,
$Y$ its dual space, $(E,\norm\cdot)$ a normed space, and $X$
a convex subset of $Y=Z'$, that is *weakly closed and unbounded or open.
Fix $K\subset X$ *weakly closed and bounded and 
$g\colon X\to\R$ a function satisfying the following conditions:

\begin{enumerate}
\item[$(g_0)$] $g$ is strictly positive, strictly convex, *weak-to-norm continuous
	on $X$  and Fr\'echet differentiable on $K$ such that
	$g'\colon K\to Y'$ is *weak-to-norm continuous and $g'(K)\subset \pi(Z)$.
\end{enumerate}
%%The mentioned before necessary growth assumptions on $g$
We make the following growth hypotheses on $g$:
\begin{enumerate}
\item[$(g_1)$] for every $n\ge1$ there exists a *weakly closed, convex and bounded subset
	$B_n$ of $X$ containing $K$ such that
	for every $t\in X\setminus B_n$ one has $g(t)\ge n$
	(or equivalently, for every $n\ge1$ setting $B_n\decl g^{-1}([0,n])$ and
	requiring that $K\subset B_n$, $B_n$ is bounded
	and $X\setminus B_n\neq \emptyset$).
	In case $X$ is unbounded, we additionally require
	\be \lim_{\stackrel{\scriptstyle \norm t\rightarrow\infty}{t\in X}}
		\frac{g(t)}{\norm{t}}=+\infty.\label{eqn:ipo.cresc.inf.g}\ee
\end{enumerate}

Define the function $h\colon K\times X\to\R$ by setting
\be h(t,u)\decl g(u)-\left[g(t)+\langle g'(t),u-t\rangle\right].\label{def:h}\ee

If the hypothesis $(g_0)$ holds,
by the *weak-to-norm continuity of $g'$ and the strict convexity of $g$,
$h$ is *weak-to-norm continuous and strictly positive for $u\neq t$.

\medskip

In the remaining part of this section we state and prove
two Korovkin-type theorems for
*weak-to-norm continuous mappings with growth prescribed by $g$.

\bt\label{teo:conv.vett.cresc} Let $Z$, $Y$, $E$, $X$, $K$, $g$ and $h$ be
	as above and there holds the conditions $(g_0)$ and $(g_1)$.
	For each $n\ge1$ let $L_n\colon D(L_n)\to\F(K;E)$
	be a linear operator dominated by a linear positive operator
	$S_n\colon D(S_n)\to\F(K)$, with $\K(X;E,g)\subset D(L_n)$ and $\K(X,g)\subset D(S_n)$.
%%	positive operator as from Definition  \ref{def:oper.dominato}.

	Then for every $t\in K$ the following statements are equivalent:
\begin{enumerate}
\item[a)]  For every $c\in E$,
	\[L_n(c)(t)\rightarrow c,\ S_n(\uni)(t)\rightarrow 1\ \mathit{and\ }
			S_n(h(t,\cdot))(t)\rightarrow 0.\]
\item[b)] For every $c\in E$, and every continuous linear
	functional $\phi\in \pi(Z)$,
	\[L_n(c)(t)\rightarrow c,\ S_n(\uni)(t)\rightarrow 1,\ 
		S_n(\phi_{|_X})(t)\rightarrow \phi(t)\ \mathit{and\ }S_n(g)(t)\rightarrow g(t).\]
\item[c)] For every $F\in\K(X;E,g)$ and $f\in\K(X,g)$,
	\[L_n(F)(t)\rightarrow F(t)\ \mathit{and\ } S_n(f)(t)\rightarrow f(t).\]
\end{enumerate}
	If the convergences in a) are uniform with respect to $t\in K$ and
	with respect to $c\in E$ then c) holds uniformly for $t\in K$.

	Moreover, if the operators $L_n$ are $S_n$-regular, then the  above conditions
	are equivalent to one of the further statements
\begin{enumerate}
\item[d)] For every $F\in\K(X;E,g)$,
	\[ L_n(F)(t)\rightarrow F(t).\]
\item[e)]  For every $f\in\K(X,g)$,
	\[ S_n(f)(t)\rightarrow f(t).\]
\item[f)] For every continuous linear functional $\phi\in \pi(Z)$,
	\[S_n(\uni)(t)\rightarrow 1,\ 
		S_n(\phi_{|_X})(t)\rightarrow \phi(t)\ \mathit{and\ }S_n(g)(t)\rightarrow g(t).\]
\end{enumerate}
\et

\boss	We remark that, if $Y$ has finite dimension $m$, then denoting by
	$(\mathrm{pr}_i)_{1\le i\le m}$ the coordinate projections on $Y$,
	the above condition  b)  reduces to the following one
\begin{enumerate}{\it
\item[b')]  for every $c\in E$, and every $i:1\dots m$,
	\[L_n(c)(t)\rightarrow c,\ S_n(\uni)(t)\rightarrow 1,\ 
		S_n(\mathrm{pr}_i)(t)\rightarrow \mathrm{pr}_i(t)\ \mathit{and\ }
		S_n(g)(t)\rightarrow g(t),\]}
\end{enumerate}
	and the convergences in a) are uniform if and only if the same holds
	true for b').
%	When $Y$ is of finite dimension the equivalence of $b)$ and $b')$ follows

This follows from the fact that $(\mathrm{pr}_i)_{1\le i\le m}$ forms a base of the
	space $Y'$.
%	Now assume that $b')$ holds uniformly.
%	Then also the convergences in $a)$ are uniform, as $g'(t)$ in
%	\refe{quattro} can be written as a finite linear combination of 
%	these same projections.
\eoss

\boss If the space $Z$ is reflexive, it is possible to simplify the hypotheses
	dropping the ''*``, substituting $\pi(Z)$ with $Y'$ and forgetting of $Z$.
	So $X$ will be a convex subset of the real reflexive  Banach space $Y$,
	that is closed and unbounded or open;
	$K\subset X$ weakly closed and bounded;
	$g\colon X\to\R$ strictly positive, strictly convex, weak-to-norm continuous
	on $X$  and Fr\'echet differentiable on $K$ such that $g'\colon K\to Y'$
	is weak-to-norm continuous and satisfying the same growth hypotheses.
\eoss

\boss Actually, as it is easy to check from the proof of the previous theorem,
	the hypothesis on $g$ may be weakened. More precisely, if we substitute
	the hypothesis $(g_0)$ with the following:
\begin{enumerate}
\item[$(g_2)$]  $g$ is strictly positive, strictly convex, Fr\'echet differentiable
	on $K$, $g'(K)\subset \pi(Z)$, $g'(K)$ is bounded in $Y'$ and the function
	$h$, defined in \refe{def:h}, is lower semicontinuous with respect to *weak topology;
\end{enumerate}
	and leave the growth hypothesis $(g_1)$,
	in the setting of the Theorem \ref{teo:conv.vett.cresc}, with further
	hypothesis that $g, h \in D(S_n)$, we obtain the implications
	$b)\Rightarrow a) \Rightarrow c)$. Moreover if the operator $L_n$ are $S_n$-regular,
	then we have the further implications
	$f)\Rightarrow b)\Rightarrow a) \Rightarrow c)\Leftrightarrow d)$.
\eoss

\bt\label{teo:stime} In the same setting of Theorem \ref{teo:conv.vett.cresc}
	assume in addition that $Y$ has finite dimension and that
	$\psi_t^2\in\C(X,g)$ for some (and hence for all) $t\in Y$.
	If $K$ is convex and $K_1\subset \stackrel\circ{K}$ is a closed subset,
	then for any $F\in\C(X;E,g)$ there exists a constant $M>0$ depending
	only on $F$,  $K$, $K_1$ and $g$ such that the estimate
	\begin{eqnarray}\norm{L_n(F)(t)-F(t)}\alli\le\alli\norm{L_n(F(t))(t)-F(t)}
		+\omega(F,\delta)(S_n(\uni)(t)\label{estimate}\\
	&& +\delta^{-2} S_n(\psi_t^2)(t))+M S_n(h(t,\cdot))(t)\nonumber\end{eqnarray}
	holds for all $\delta>0$ and $t\in K_1$ (here
	$\omega(F,\cdot)$ stands for the modulus  of continuity of $F$ on $K$).
	When $L_n$ preserves the constants and $S_n(\uni)(t)=1$, the above
	 estimate becomes:
	\be \norm{L_n(F)(t)-F(t)}\le 2\omega(F,\gamma_n(t))+MS_n(h(t,\cdot))(t).
			\label{estimate2}\ee
	Finally, if $S_n$ preserves the linear functionals, then
	\be\norm{L_n(F)(t)-F(t)}\le 2\omega(F,\gamma_n(t))+M(S_n(g)(t)-g(t)).\label{ultima}\ee

	In case $\dim(Y)=1$, $X=\lb a,+\infty\lb$ \risp{ X=$\rb -\infty,b\rb$}
	and $K=\lb a,b\rb$, the previous estimates hold with
	$K_1=\lb a,b_1\rb$ for any $b_1<b$ \risp{ $K_1=\lb a_1,b\rb$ with
	$a<a_1$}.
\et

Before proving the theorems, we present two useful lemmas:

\bl\label{lem:stime} Let $Z$, $Y$, $E$, $X$, $K$, $g$ and $h$ be as in the
	Theorem \ref{teo:conv.vett.cresc} and consider $F\in\K(X;E,g)$.
	Then there exist an integer $\nu\ge 1$ and a constant $M>0$ such that
	\be \norm{F(t)-F(u)}\le Mh(t,u)\quad \mathrm{for\ any\ }t\in K\ \mathrm{and\ }u\in X\setminus B_\nu.\label{uno}\ee

	Moreover, for any $\delta>0$ and any finite set $\ell\subset Z$
	one gets
	\be \norm{F(t)-F(u)}\le \omega(F,K,I_{\ell,\delta})+Mh(t,u)
	\quad \mathrm{for\ any\ }t\in K\ \mathrm{and\ }u\in X,\label{due}\ee
	where $I_{\ell,\delta}$ is the following neighborhood of $0$ in the
	 *weak topology on $Y$:
	\[I_{\ell,\delta}\decl\{y\in Y\big|\forall \xi\in\ell: |y(\xi)|<\delta\},\]
	and
	\be\omega(F,K,I_{\ell,\delta})\decl\sup\{\norm{F(t)-F(u)}\,|t\in K,\,
		u\in B_\nu,\,u\in  t+I_{\ell,\delta}\}\label{def:omega}\ee
\el
\bp \textbf{The estimate \refe{uno}.} From the *weak-to-norm continuity of $F$, $g$,
	$g'$ and the boundedness of $K$, it follows that there exists
	a positive constant $M_1>0$
	such that for all $t\in K$ one has $\norm{F(t)}\le M_1$, $|g(t)|\le M_1$,
	$\norm{g'(t)}_{Y'}\le M_1$ and $\norm t\le M_1$.
	Thus for $t\in K$ and $u\in X$ we get
	\[\frac{\langle g'(t),u-t\rangle}{g(u)}\le\frac{M_1\norm{u-t}}{g(u)}\le
		 \frac{M_1}{g(u)}(\norm u+M_1)\]
	and then
	\[\frac{h(t,u)}{g(u)}\ge1-M_1\frac{1+\norm u+M_1}{g(u)}\quad
			\mathrm{for\ all\ } t\in K\ \mathrm{and\ }u\in X.\]
	Hence, by the hypotheses on the growth of $g$, it follows that
	\be 0<M_1\frac{1+\norm u+M_1}{g(u)}\le\frac{M_1+M_1^2}{n}+\frac{\norm u}{g(u)}\quad
			\mathrm{for\ any\ } u\in X\setminus B_n.\label{tec}\ee
	Fix $\epsilon\in\rb 0,1\lb$. If $X$ is bounded, that is $\norm u \le N$
	for $u\in X$ and some constant $N$, then taking $n$ greater than
	an appropriate integer $\nu$ we obtain
	\[M_1\frac{1+\norm u+M_1}{g(u)}\le M_1\frac{1+N+M_1}{g(u)}\le
		M_1\frac{1+N+M_1}{\nu}<\epsilon\]
	for all $n\ge\nu$ and $u\in X\setminus B_n$.
	If $X$ is unbounded, then by \refe{eqn:ipo.cresc.inf.g}, there exists
	$a>0$ such that for any $\norm u\ge a$ we have
	$\frac{\norm u}{g(u)}<\epsilon/2$.
	Setting  $\nu\decl 2 \max\{a,M_1+M_1^2\}/\epsilon$, for any
	$n\ge\nu$ and $u\in X\setminus B_n$, one has
	$\frac{\norm{u}}{g(u)}\le\epsilon/2$
	(in both cases $\norm{u}\ge a$ and $\norm{u}<a$), then looking at
	(\ref{tec}) we obtain
%	corrispondenza di $a$, posto
	\[ M_1\frac{1+\norm u+M_1}{g(u)}\le\frac{M_1+M_1^2}{n}+\frac{\norm u}{g(u)}\le
		\frac \epsilon 2 + \frac \epsilon 2= \epsilon \]
	as in the case of $X$ bounded.
	Hence for $n\ge\nu$,
%%	$1-M_1\frac{1+\norm u+M_1}{g(u)}$ does not vanish and for
	$u\in X\setminus B_n$ and $t\in K$ we have
	\bern\frac{\norm{F(t)-F(u)}}{h(t,u)}&=&\frac{\norm{F(t)-F(u)}}{g(u)}\frac{g(u)}{h(t,u)}\\
		&\le&\frac{\norm{F(u)}+M_1}{g(u)}\left(1-M_1\frac{1+\norm u+M_1}{g(u)}\right)^{-1}\\
		&\le&\frac{\norm{F(u)}+M_1}{g(u)}(1-\epsilon)^{-1}.\eern
	The above inequality together with $\norm{F(u)}\le Mg(u)$ accomplishes
	the proof of \refe{uno}.

\medskip
\noindent\textbf{The estimate \refe{due}.} Set
	\[ A\decl\{(t,u)\big|t\in K,\,u\in B_\nu\, \mathrm{and\ } u\not\in t+I_{\ell,\delta} \}.\]
	$A$ is *weakly closed and bounded, because the same holds for $K$ and
	$B_\nu$.
	Since $h$ is *weak-to-norm continuous, then by Weierstrass' theorem, we deduce
	that $h$ has a minimum $m$ on $A$, and $m>0$ because $h(t,u)=0$ only for
	$u=t$. Moreover, since $F$ is *weak-to-norm continuous, the same holds
	true for the function
	$\norm F$, and, consequently, $\norm F$ is bounded on the bounded set $B_\nu$.
	Hence we obtain
	\[ \norm{F(t)-F(u)}\le 2\norm{F}_{B_\nu} \frac{h(t,u)}{m}=M_2h(t,u)\]
	for every $t\in K$ and $u\in B_\nu\setminus (t+I_{\ell,\delta})$.

	Recalling the estimate \refe{uno} and the definition \refe{def:omega},
	we conclude the proof of \refe{due}.
\ep

The next lemma explains an important property of $\omega(F,K,I_{\ell,\delta})$,
that will be used in the sequel.

\bl Under the some assumptions of Lemma \ref{lem:stime}, it follows that
	for any positive real $\epsilon>0$ there exist a
	finite set $\ell\subset Z$ and a constant $\delta>0$ such that
	$\omega(F,K,I_{\ell,\delta})\le\epsilon$.
\el
\bp	By the *weak-to-norm continuity of $F$, for a fixed $t\in K$ there exist
	a finite set $\ell_t\subset Z$ and $\delta_t>0$ such that
	$\norm{F(t)-F(u)}<\epsilon/2$ for $u\in t+I_{\ell_t,\delta_t}$.
	Trivially $K\subset\bigcup_{t\in K} t+I_{\ell_t,\delta_t/2}$.
	Since $K$ is compact in the *weak topology,
	there are $t_1,t_2,\dots,t_n \in K$, such that
	\[K\subset\bigcup_{i=1}^n t_i+I_{\ell_i,\delta_i/2}.\]
%	where, for short,
%	we set $\delta_i=\delta_t$ ed $\ell_i=\ell_t$ for $t=t_i$, $i=1\dots n$.
	Let $\delta\decl 1/2\min\{\delta_i,\,i=1\dots n\}$ and 
	$\ell\decl\bigcup_{i=1}^n\ell_i$. We prove that $I_{\ell,\delta}$ is
	the desired neighborhood of zero.

	Fix  $t\in K$ and $u\in (t+I_{\ell,\delta})\cap B_\nu$.
	Let $i$ be the index for which $t\in t_i+I_{\ell_i,\delta_i/2}$.
	For any $\xi\in\ell_i$ the inequality
	\[|\xi(u-t_i)|\le|\xi(u-t)|+|\xi(t-t_i)|<\delta+\delta_i/2\le\delta_i\]
	holds, and thus $u\in t_i+I_{\ell_i,\delta_i}$.
	Therefore
	\[\norm{F(t)-F(u)}\le\norm{F(t)-F(t_i)} + \norm{F(t_i)-F(u)}<\epsilon/2
	+\epsilon/2,\]
	which yields precisely the desired estimate for
	$\omega(F,K,I_{\ell,\delta})$.
\ep

\medskip
Now we prove our main results.
\medskip

\noindent \textbf{Proof of Theorem \ref{teo:conv.vett.cresc}.\ } First of all,
	observe that for given $F\in\K(X;E,g)$ and $t\in K$,
	applying $S_n$  to both sides of \refe{due} of Lemma \ref{lem:stime},
	we obtain
	\[ S_n(\norm{F-F(t)})(t)\le S_n(\uni)(t)\omega(F,K,I_{\ell,\delta})
		+MS_n(h(t,\cdot))(t).\]
	Consequently
	\begin{eqnarray}\norm{L_n(F)(t)-F(t)}&\le&\norm{L_n(F(t))(t)-F(t)}+
				S_n(\norm{F-F(t)})(t)\label{tre}\\
	&\le&\norm{L_n(F(t))(t)-F(t)}+S_n(\uni)(t)\omega(F,K,I_{\ell,\delta})\nonumber\\
	&&\qquad +MS_n(h(t,\cdot))(t).\nonumber\end{eqnarray}

	We prove the implication a)$\Rightarrow$ c).
%	In order to obtain the convergence of $L_n(F)(t)$ to $F(t)$, we
	Take $\epsilon>0$ %% In correspondence of $\epsilon/6$ we
	and consider the zero neighborhood
	$I_{\ell,\delta}$ for which $\omega(F,K,I_{\ell,\delta})\le\epsilon/6$.
	By Lemma \ref{lem:stime},
	there exists a constant $M$ such that the relation \refe{due} holds
	for  $I_{\ell,\delta}$.
%	For the hypothesis on the operators $L_n$ and $S_n$ let
	In view of a), for $n$ sufficiently large we have
	$S_n(h(t,\cdot))(t)<\epsilon/(3M)$,
	$S_n(\uni)(t)<2$ and $\norm{L_n(F(t))(t)-F(t)}<\epsilon/3$,
	and thus, using \refe{tre} we deduce
	\[\norm{L_n(F)(t)-F(t)}\le\epsilon/3+2\epsilon/6+M\epsilon/(3M)=\epsilon,\]
	that proves the convergence of $L_n(F)(t)$ to $F(t)$.
	It is clear that the convergence is uniform if the same holds for a).

	Fix $f\in\K(X,g)$. In order to prove the convergence of  $S_n(f)(t)$
	to $f(t)$ we proceed in the manner we made before
	substituting the norm $\norm{\cdot}$ in $E$  with the absolute value.
% Dall'osservazione che la forte continuit\`a e la debole
%	continuit\`a per le funzioni di $\C(X;\R)$ sono la stessa cosa,
%	Possiamo rifare per $f$ gli stessi ragionamenti, gi\`a fatti per $F$,
%	a patto di sostituire la norma di $E$ con il valore assoluto e
%	and to have the convergence of (uniform if need be).

\bigskip
	In order to prove the implication $c)\Rightarrow b)$, it is sufficient
	to observe that the constant functions are *weak-to-norm continuous,
	and  the function $g$ and all continuous functionals in $\pi(Z)$
	belongs to $\K(X,g)$ (by \refe{eqn:ipo.cresc.inf.g}).

% Le prime due convergenze sono vere per ipotesi.
%	Per provare le altre due basta applicare l'implicazione
%	$a)\Rightarrow c)$ gi\`a provata agli operatori $S_n$, essendo questi dominati
%	banalmente da se stessi ed osservando che $g\in\C(X,g)$ e che se $\phi_{|_X}$
%	\`e la restrizione ad $X$ di una forma lineare continua, per l'ipotesi
%	\refe{eqn:ipo.cresc.inf.g}, $\phi_{|_X}\in\C(X,g)$. 
\bigskip
	The implication $b)\Rightarrow a)$ follows directly from the identity
	\be S_n(h(t,\cdot))(t)=S_n(g)(t)-g(t)S_n(\uni)(t)-S_n(g'(t))(t)+
			\langle g'(t),t\rangle S_n(\uni)(t).\label{quattro}\ee

	Now we assume that $L_n$ is $S_n$-regular.
	The implication  $d)\Rightarrow c)$.
	Fix $f\in \K(X,g)$. Taking $x\in E$, by definition of $S$-regularity,
	we have
	\[ S_n(f)(t)\otimes x =L_n(f\otimes x)(t),\]
	that converges to $f(t)x$. Since $x$ is arbitrary we have the convergence
	of $S_n(f)(t)$ to $f(t)$.
%%%%%%	e) implies b)
%	let $c\colon X\to E$ be an constant map of value $c$. We have
%	\[\norm{L_n(c)(t)-c}=\norm{\int_X c\de_u\mu_n(u;t)-c}
%			=\norm c\abs{S_n(\uni)(t)-1}\rightarrow 0,\]
%	so the first convergence of the proposition b) holds. The other ones
%	are true by hypotheses, because $\uni,\ \phi_X$ and $g$ belongs to $\K(X,g)$.
%\bigskip
%	When $Y$ is of finite dimension the equivalence of $b)$ and $b')$ follows
%	from the fact that $(\mathrm{pr}_i)_{1\le i\le m}$ forms a base of the
%	space $Y'$.
%	Now assume that $b')$ holds uniformly.
%	Then also the convergences in $a)$ are uniform, as $g'(t)$ in
%	\refe{quattro} can be written as a finite linear combination of 
%	these same projections.
	The implication $f)\Rightarrow b)$ follows from identity
	\[L_n(c)(t)=L_n(\uni\otimes c)(t)=S_n(\uni)\otimes c \]
	and the missing implication $e)\Rightarrow f)$ is immediate. The proof is complete.
{\hfill $\Box$}

\bigskip

\noindent \textbf{Proof of Theorem \ref{teo:stime}.\ } Fix $F\in\K(X;E,g)$
	and $\delta>0$.
	By \refe{uno} for every $t\in K_1$ and $u\in X\setminus B_\nu$ we get
	\be \norm{F(t)-F(u)}\le M_1h(t,u).\label{cinque}\ee

	On the other hand the inequality 
	\be \norm{F(t)-F(u)}\le \omega(F,\norm{t-u})\le
			(1+\delta^{-2} \norm{t-u}^2)\omega(F,\delta)\label{sei}\ee
	holds for every $t\in K_1$ and $u\in K$ ($\omega(F,\delta)$ stands here
	for the modulus of continuity of $F$ on $K$).

	Now we discuss the case $t\in K_1$ and $u\in \overline{B_\nu\setminus K}$.
	Since $K_1\subset\stackrel{\circ}K$, there exists a closed and convex
	set $K_\eta\subset\stackrel{\circ}K$ such that
	$K_1\subset\stackrel{\circ}{K_\eta}$.
	From the convexity of $B_\nu$, $K$ and $K_\eta$, it follows that
	\[[a',a'']=[u,t]\cap\overline{K\setminus K_\eta}\]
	for some $a'\in K$ and $a''\in K_\eta$.
	Let $P\colon\lb 0,1\rb\to[u,t]$ be the parametric representation
	of the segment, $P(s)\decl(1-s)u+st$ ($0\le s\le1$),
%	 $P(0)=u$, $P(1)=t$ 
	and  $0\le s'\le s''\le1$ such that $P(s')=a'$ and $P(s'')=a''$.
	We set $\hat g\decl g\circ P\colon\lb0,1\rb\to[u,t]$ and
	$\hat h(r,s)\decl\hat g(s)-[\hat g(r)+\hat g'(r)(s-r)]$.
	Note that $\hat g$ is strictly convex by the strict convexity of $g$.
	This yields $\hat h(s'',s')\le\hat h(1,0)$.
	Observing that
	\[\hat g'(r)=\langle g'(P(r)),P'(r)\rangle=\langle g'(P(r)),t-u\rangle,\]
	and $P(s)-P(r)=(s-r)(t-u)$, we get: 
	\[\hat h(r,s)=g(P(s))-[g(P(r))+\langle g'(P(r)),P(s)-P(r)\rangle]=h(P(r),P(s)).\]
	Hence
	\[h(a'',a')=\hat h(s'',s')\le\hat h(1,0)= h(t,u),\]
	and consequently
	\be \norm{F(t)-F(u)}\le \norm{F(t)}+\norm{F(u)}\le2\norm{F}_{B_\nu} 
		h(t,u)/ h(a'',a').\label{sette}\ee
	Since $\partial K\cap\partial K_\eta=\emptyset$, surely
	$\inf\{h(a'',a')| a'\in\partial K,\, a''\in\partial K_\eta \}>0$ and therefore
	\be\norm{F(t)-F(u)}\le M_2h(t,u)\quad\mathrm{for\ any\ }t\in K_1,\,u\in \overline{B_\nu\setminus K}.\label{otto}\ee

	In case $\dim(Y)=1$, $X=\lb a,+\infty\lb$, $K=\lb a,b\rb$ and
	$K_1=\lb a,b_1\rb$ (with $b_1<b$) the relation \refe{otto} is
	established in the similar manner.
	One considers $K_\eta\decl\lb a,b_2\rb$ with $b_1<b_2<b$ and finds
	$a'=b$ and $a''=b_2$, which yield \refe{otto} in view \refe{sette}
	and the inequality $0<h(b_2,b)\le h(t,u)$.

	Combining inequalities \refe{cinque}, \refe{sei} e \refe{otto} we
	obtain
	\[\norm{F(t)-F(u)}\le(1+\delta^{-2} \norm{t-u}^2)\omega(F,\delta)+Mh(t,u)\]
	for all $t\in K_1$ and $u\in X$.
	Now applying $S_n$ and using the first inequality in \refe{tre}
	we obtain estimate \refe{estimate}.
	The last inequality \refe{ultima} easily follows from relation \refe{quattro}.
{\hfill $\Box$\par\medskip}

\boss We stress the fact that the constant $M$ in \refe{estimate},
	\refe{estimate2} and \refe{ultima} depends only on
	$F$, $K$, $K_1$ and $g$; in particular, it does
	\emph{not} depend on the operators $L_n$ or $S_n$.
\eoss
%% Questo ci permetter\`a
%%nel prossimo paragrafo di ottenere ulteriori processi di approssimazione
%%con stime uniformi sui compatti.

\boss From the previous theorems we deduce that an approximation process for
	real valued functions $S_n$, defined by means of positive measures,
	yields another process $L_n$, for vector valued functions.
	Note that the process $L_n$ ``inherits'' the estimates valid for $S_n$.
\eoss

Theorem \ref{teo:conv.vett.cresc} and \ref{teo:stime} generalize the
corresponding results in \cite{shaw} and \cite{shaw-yeh}.
The main result in \cite{ditzian} is an easy consequence of our
Theorem \ref{teo:stime} under the additional requirement that
the control function $(1+t^2)\mu(t)$ is strictly convex.

Moreover, Theorem \ref{teo:stime} extends the results of \cite{shaw},
providing them with estimates of the corresponding rate of convergence.